\documentclass[11pt]{article}

\usepackage{amssymb}
\usepackage{amsthm}
\usepackage[english]{babel}
\usepackage{mathtools}
\usepackage{graphicx}
\usepackage[all,cmtip]{xy}
\usepackage{tikz}
\usetikzlibrary{matrix,arrows}
\usepackage{float}

\theoremstyle{plain}
\newtheorem{theorem}{Theorem}
\newtheorem{lemma}[theorem]{Lemma}
\newtheorem{corollary}[theorem]{Corollary}

\theoremstyle{definition}
\newtheorem{definition}[theorem]{Definition}

\theoremstyle{remark}

\title{ Quasi-projective monounary algebras}

\author{\' Eva Jung\'abel}

\date{}

\begin{document}

\maketitle

\begin{abstract}
Wu and Jans introduced quasi-projective modules where they say a $\cal R$ module $\cal M$ is quasi-projective if  for every submodule $\cal N$, for every homomorphism $f : {\cal M} \rightarrow {\cal M}/{\cal N}$ and every epimorphism $j: {\cal M}\rightarrow {\cal M}/{\cal N}$ there is an endomorphism $\phi$ of $\cal M$ such that $\phi\circ j=f$. We say that a structure $\cal S$ is quasi-projective if for every structure $\cal T$, for every homomorphism $f : {\cal S} \rightarrow {\cal T}$ and every epimorphism $j: {\cal S}\rightarrow {\cal T}$ there is an endomorphism $\phi$ of $\cal S$ such that $\phi\circ j=f$. In 2004 D. Jakub\'ikov\'a-Studenovsk\'a defined the concept of the factor algebra denoted by ${\cal A}/{\cal B}$, where ${\cal A}$ is a monounary algebra and ${\cal B}$ is a subalgebra of $\cal A$. In this paper, we characterise the quasi-projective monounary algebras of arbitrary cardinalities for the definition of D. Jakub\'ikov\'a-Studenovsk\'a and for the second definition.

\end{abstract}

\let\thefootnote\relax

\section{Introduction}\label{intro}

By a structure we mean a set together with an indexed set of relations and operations on it. A first-order structure $\cal A$ is called {\it homogeneous} if any isomorphism between two finitely generated substructures of $\cal A$ is induced by some automorphism of $\cal A$. In several classes of combinatorial structures the homogeneous structures are classified. 

P. Cameron and J. Ne\v set\v ril \cite{cameron-nesetril} introduced the following variant of homogeneity: a structure is called {\it homomorphism-homogeneous} if every homomorphism between finite induced substructures extends to an endomorphism of the structure. 

Ho\-mo\-mor\-phism-ho\-mo\-ge\-neous graphs were investigated by Rusinov and Schweitzer  in \cite{rusinov-schweitzer}. Among others it is shown that the problem of deciding if the graph is homomorhism-homogeneous is coNP-complete. Finite algebras also harbour some classes of high computational complexity \cite{masulovic3}, hence we cannot expect a brief classification in case of algebraic structures in general. A characterisation of all ho\-mo\-mor\-phism-ho\-mo\-ge\-neous partial orders of arbitrary cardinalities with non-strict relation is given  by Ma\v sulovi\' c  \cite{masulovic2} and {Cameron and Lockett} \cite{cameron-lockett}, independently. Several  other ho\-mo\-mor\-phism-ho\-mo\-ge\-neous structures are characterised including monounary algebras of arbitrary cardinalities by Ma\v sulovi\'c and Jung\'abel \cite{jungabel-masulovic}.

\vskip3mm
Interestingly, the concept of homomorphism-homogeneity is not a recent concept, it exists under the name of {\it quasi-injectivity} just with the slight difference. A structure is said to be quasi-injective if every homomorphism from an arbitrary substructure of the structure into the structure extends to an endomorphism of the structure.

All quasi-injective Abelian groups are described \cite{fuchs1} as finite quasi-injective groups \cite{bertholf-walls}. Infinite quasi-injective groups \cite{tomkinson} are partly characterised. There are results about quasi-injective modules by Johnson and Wong \cite{johnson-wong1}, Harada \cite{harada1}, Faith and Utumi \cite{faith-utumi1}, Fuchs \cite{fuchs2} and others.

The dual concept of quasi-injectivity is {\it quasi-projectivity} and it was introduced for modules by Wu and Jans \cite{wu-jans1} in 1967. Phrased in terms of diagrams, the module $\cal S$ is quasi-injective if every diagram

\centerline{
\xymatrix@=5em{
  0 \ar[r] 
& {\cal T} \ar[r]^j \ar[d]^f
&{\cal S}\\
& {\cal S}}
}

\noindent
can be embedded in a commutative diagram

\centerline{
\xymatrix@=5em{
  0 \ar[r] 
& {\cal T} \ar[r]^j \ar[d]^f
&{\cal S} \ar@{-->}[dl]_{\phi} \\
& {\cal S}}
}

\noindent
where $\cal T$ is a submodule of $\cal S$, $j$ is the monomorphism  and $f$ is a homomorphism of $\cal T$ into $\cal S$.

The module $\cal S$ is said to be quasi-projective if every diagram

\centerline{
\xymatrix@=5em{
& {\cal S} \ar[d]^f \\
{\cal S} \ar[r]^j
& {\cal S/T} \ar[r]
& 0}
}

\noindent
can be embedded in a commutative diagram

\centerline{
\xymatrix@=5em{
& {\cal S} \ar[d]^f \ar@{-->}[dl]_-{\phi} \\
{\cal S} \ar[r]^j
& {\cal S/T} \ar[r]
& 0}
}

\noindent
where $\cal T$ is a submodule of $\cal S$, $j$ is the epimorphism  and $f$ is a homomorphism of $\cal T$ into $\cal S$.

In \cite{wu-jans1} some properties of quasi-projective modules are shown, a structure theorem for indecomposable finitely generated quasi-projectives over semi-perfect rings is obtained and the finitely generated quasi-projective Abelian groups are described. The general case of quasi-projective Abelian groups is characterised by Fuchs and Rangaswamy \cite{fuchs-rangaswamy1}. A decomposition theorem that is a characterisation for quasi-projective modules over left perfect rings is given by Koehler \cite{koehler}. 

We say that a structure $\cal S$ is quasi-projective if for every structure $\cal T$, for every homomorphism $f : {\cal S} \rightarrow {\cal T}$ and every epimorphism $j: {\cal S}\rightarrow {\cal T}$ there is an endomorphism $\phi$ of $\cal S$ such that $\phi\circ j=f$. D. Jakub\'ikov\'a-Studenovsk\'a \cite{jakubikova1} in 2004 defined the concept of the factor algebra denoted by ${\cal S}/{\cal T}$, where ${\cal S}$ is a monounary algebra and ${\cal T}$ is a subalgebra of $\cal S$. In this work we characterise the quasi-projective monounary algebras of arbitrary cardinalities for the definition of D. Jakub\'ikov\'a-Studenovsk\'a and for the second definition.

\section{Preliminaries}\label{preli}

A monounary algebra (see \cite{jakubikova-pocs}) is an algebra ${\cal A}=(A,\alpha)$ where $\alpha \colon A \rightarrow A$ is an unary operation on $A$. For $\emptyset \neq T \subseteq A$, the {\it subalgebra} of $(A,\alpha)$ generated by $T$ is the algebra $(\langle T \rangle, \alpha|_{\langle T \rangle})$, where $\langle T \rangle =\{ \alpha^k(t) \mid k \geqslant  0, t \in T \}$. 

We define a binary relation $\sim$ on $A$ as follows: $a \sim b$ if there exists a number $k$ such that $\alpha^k(a)=b$ or $\alpha^k(b)=a$. It is easy to show that $\sim$ is an equivalence relation and elements of $A/ _\sim$  are referred to as {\it connected components} of $(A,\alpha)$.

An element $a \in A$ is {\it cyclic} if there exists a $k \geqslant 1$ such that $\alpha^k(a) = a$. Otherwise, $a$ is said to be {\it acyclic}. The set of all cyclic elements in a connected component $S \subseteq A$ is called {\it the cycle of} $S$. It may happen that a connected component does not have a cycle.  The {\it length of a cycle} $C$ is the least number $k$ such that $\alpha^k(c)=c$ for all $c\in C$.

For a connected component $S\subseteq T$ let $cn(S)$, the cycle number of $S$, denote the length of the cycle in $S$. If $S$ does not have a cycle, we set cn~$(S)=\infty$. 

Let $C$ be a cycle. The {\it distance}, $d(c_i,c_j)$, between two elements $c_i, c_j\in C$, where $i<j$, is the least number $k$ such that $\alpha^k(c_i)=c_j$. Note that the function $d$ is not symmetric. 

An element $a\in A$ is called a {\it leaf} in $(A,\alpha)$ if $\alpha^{-1} (a) = \emptyset$.  Let $I$ be one of the sets $\{1,2,\ldots, n\}$,  ${\mathbb N}$, $\mathbb{Z}^- \cup \{0, 1,2,\ldots, n\}$ or ${\mathbb Z}$. A {\it branch $B=\{b_i\}_{i\in I}$} in a monounary algebra $\cal A$ is a maximal sequence of $b_i$, $i\in I$ such that $b_i = \alpha(b_{i+1})$ for all $i \in I$ and $b_i$ is acyclic for all $i \in I$. A branch can be finite or infinite. If a branch is finite, then it has a leaf. A branch can end at a cycle or not. We say that a branch  $B=\{b_1, b_2, b_3, \ldots\}$ ends at a cycle $C$ if $\alpha(b_1)\in C$ and we say that $\alpha(b_1)$ is {\it the ending point } of the branch. 

Let $b$ be an acyclic element, then $\uparrow b$, the {\it bunch} of $b$, denote the set $\{\alpha^{-i}(b)\}_{i\geq 0}$. Let $B$ be a branch that ends at a cycle $C$. Then $\uparrow{ B}$, the {\it bunch} of $B$, denote the set $\{\alpha^{-i}(b_1)\}_{i\geq 0}$, where $b_1\in B$ and $\alpha(b_1)\in C$.

For an acyclic element $a \in A$, let $l_h(a)$, the {\it height of a}, denote the least $k \geqslant 1$ such that $\alpha^k(a)$ is a cyclic element. If no such $k$ exists, we set $l_h(a) = \infty$. For a branch $B$ with a leaf $a$, $l_h(B)=l_h(a)$. 

For an acyclic element $a \in A$ let $l_d(a)$, the {\it depth of a}, denote the greatest $k \geqslant 0$ such that $\alpha^{-k}(a)$ is a leaf. If no such $k$ exists, we set $l_d(a) = \infty$. For a branch $B$ and for an element $b_g\in B$, $l^g_d(B)$, the {\it depth of $B$}, denote the number $k \geqslant 0$ such that $\alpha^{-k}(b_g)$ is a leaf and $\alpha^{-i}(b_g)\in B$, for all $i\in \{0,1,\ldots, k\}$. If no such $k$ exists, we set $l^g_d(B) = \infty$. 

For finite branch which ends in a cycle, we have $l_h(B)=l_d(B)$. 
Let $d(b_1,b_2)$, the {\it distance} between two elements $b_1$ and $b_2$ in a branch $B$, denote the unique number $k$ such that $\alpha^k(b_1)=b_2$ or $\alpha^k(b_2)=b_1$. 

The diagram of a monounary algebra is a graph with loops where vertices are the elements and an edge is between elements $a_i$ and $a_j$  if $\alpha(a_i)=a_j$.

A mapping $f \colon {\cal A}_1 \rightarrow {\cal A}_2$ is a {\it homomorphism} from $({\cal A}_1,\alpha_1)$ to $({\cal A}_2,\alpha_2)$ if $\alpha_2 \circ f =f \circ \alpha_1$.

\begin{lemma}
Let $f$ be an endomomorphism of a monounary algebra $(A,\alpha)$. Then:

\begin{enumerate}
\item every cyclic element is mapped to a cyclic element and the cycle $C_1$ is mapped to the cycle $C_2$  if and only if cn~$(C_2)|$cn~$(C_1)$. If the there is a branch $B_1$ which ends in the cycle $C_1$, then it can be mapped to a cycle $C_2$ or to a branch $B_2$ which ends at a cycle $C_2$,

\item a branch which does not end at a cycle and which has a leaf can be mapped to cycle with branches with leaves or without leaves or to a branch  without cycle and with or without leaf by homomorphism. Also, a branch which does not end at a cycle and which has not a leaf can be mapped to cycle with branches with leaves or without leaves or to a branch without cycle and without leaf by homomorphism.

\end{enumerate}
\end{lemma}

\vskip3mm
According to D. Jakub\'ikov\'a-Studenovsk\'a, for a subalgebra $\cal U$ of al monounary algebra $\cal A$ we define the quotient
monounary algebra ${\cal A}/{\cal U}$.

For an equivalence relation $\theta$ on ${\cal A}=(A, \alpha)$, it is called a congruence of ${\cal A}$ if $x, y \in A$, $(x, y) \in \theta$ implies $(\alpha(x), \alpha(y)) \in \theta$. For $x \in A$ , the equivalence class of $\theta$ containing $x$ is denoted by $[x]_{\theta}$ or $[x]$. A quotient algebra ${\cal A}_{\theta} = (A_{\theta}, \alpha_{\theta})$ is such that $A_\theta$ is the union of equivalence classes and $\alpha_{\theta}([x]) = [\alpha(x)]$.

\begin{definition}[D. Jakub\'ikov\'a-Studenovsk\'a \cite{jakubikova1}]
Let ${\cal A} = (A, \alpha)$, $\emptyset \not= U \subseteq A$. We denote by $\theta_U$ the smallest congruence relation of $\cal A$ such that if $x, y \in U$ belong to the same connected component of $\cal A$ , then $x, y$ belong to the same equivalence class of the congruence $\theta_U$.
\end{definition}

\begin{lemma}[D. Jakub\'ikov\'a-Studenovsk\'a \cite{jakubikova1}]
Suppose that $\cal A = (A,\alpha)$ , ${\cal U} = (U, \theta_U)$ is a subalgebra of ${\cal A}$ . Let $x, y \in A$. Then $(x, y) \in \theta_U$ if and only if either $x, y$ belong to the same connected component of $\cal A$ and $\{x, y\} \subseteq U$ or $x = y$.
\end{lemma}

\begin{corollary}[D. Jakub\'ikov\'a-Studenovsk\'a \cite{jakubikova1}] Let $\cal A$ be connected, and ${\cal U}=(U,\alpha_U)$ be a subalgebra of $\cal A$ , $|U| > 1$. Then the unique nontrivial equivalence class of $\theta_U$ is equal to $U$.
\end{corollary}

\begin{definition}[D. Jakub\'ikov\'a-Studenovsk\'a \cite{jakubikova1}] Let ${\cal A} = (A, \alpha)$ and let ${\cal U} = (U, \alpha_U)$ be a subalgebra of ${\cal A}$. By a {\it quotient monounary} algebra ${\cal A}/{\cal U} = (A/U, \alpha_{A/U})$ we understand algebra $A/\theta_U$.
\end{definition}

\begin{corollary}[D. Jakub\'ikov\'a-Studenovsk\'a \cite{jakubikova1}]
Let ${\cal A} = (A, \alpha)$ be connected and complete, and ${\cal U} = (U, \alpha_U)$ be its subalgebra. Then
\begin{enumerate}
\item[(i)] $\alpha_{A/U}(\{x\})= \{\alpha(x)\}$ if $x \in A$, $\alpha(x)\notin U$,
\item[(ii)] $\alpha_{A/U}(\{x\})= U$ if $x \in A$, $\alpha(x)\in U$,
\item[(iii)] $\alpha_{A/U}(U)= U$.
\end{enumerate}
\end{corollary}

\section{Monounary alegbras with a cycle}\label{qpma1}

In this subsection we characterise quasi-projective monounary algebras with a cycle. The characterisation is the following:

\begin{theorem}\label{mono2.11}
Let ${\cal A}=(A,\alpha)$ be a monounary algebra with a cycle. It is quasi-projective if and only if (Fig. \ref{Tetel41}):

\begin{enumerate}

\item\label{item11} there are cycles in all connected components,

\item\label{item12} cycles have the same length  in all connected components,

\item\label{item13} in every bunch there is one branch,

\item\label{item14} all branches have the same length,

\item\label{item15} 
\begin{enumerate}
\item[(i)] in every connected component there is at most one branch or 

\item[(ii)] there are two branches in a connected component and their cyclic element at which they end is on the same distance from each other.
\end{enumerate}
\end{enumerate}

The statement holds if we use the notion of factor algebra of Jakub\'ikov\'a-Studenovsk\'a's definition.

\end{theorem}

\begin{figure}[H]
\begin{center}
\includegraphics[width= 14 cm ]{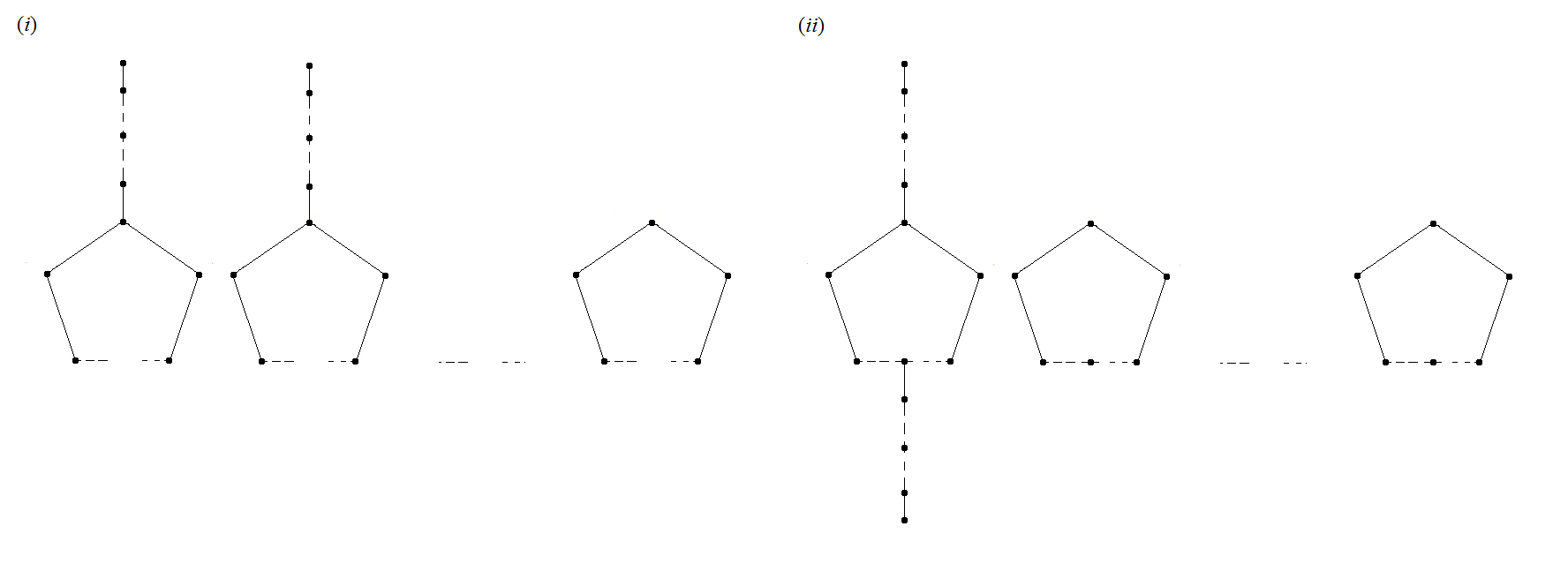}
\caption{\label{Tetel41} Theorem \ref{mono2.11}}
\end{center}
\end{figure}

Now, we start to prove Theorem \ref{mono2.11} step by step.

\begin{lemma}\label{mono2.1}
If a monounary algebra ${\cal A}=(A,\alpha)$ is quasi-projective, then there does not exist a cycle and a branch without a cycle. Moreover, if there are cycles, all cycles have the same length.

The statement and the proof holds if we use the notion of factor algebra of Jakub\'ikov\'a-Studenovsk\'a's definition choosing $U=A$.

\end{lemma}

\begin{proof}
Let ${\cal A}=(A,\alpha)$ be a quasi-projective monounary algebra. Let $S_1$ be a connected components with a cycle $C_1=\{c^1_1, c^1_2,\ldots, c^1_m\}$ in it. Suppose to the contrary that in an other connected component $S_2$ there is a branch without cycle or a cycle with a different length as $C_1$. In any case we have that cn~$(S_1)\not=$cn~$(S_2)$. We may assume that $cn(S_1)< cn(S_2)$. Let ${\cal T}=(T, \alpha_T)$ be a monounary algebra such that $T=C^T_1 \cup C^T_2$, where $C^T_1$ and $C^T_2$ are one point cycles with elements $c^T_1$ and $c^T_2$, respectively. We define $f$ and $j$ in the following way:

$$f:\left(\begin{array}{rrrrr}  S_1  & S_2 & S_i \\  C^T_1 & C^T_2  & C^T_2 \end{array} \right), i\not= 1,2,$$

$$j:\left(\begin{array}{rrrrr}   S_1  & S_2 & S_i \\    C^T_2 & C^T_1 & C^T_2 \end{array} \right), i\not= 1,2.$$

\noindent
$S_1$ is mapped to the cycle $C^T_1$ by $f$ in the following way: $f(c^1_i)=c^T_1$ for $i\in \{1,2,\ldots, m\}$. Let $B_1$, $ B_2 $, $\ldots$, $\ B_p $ be branches ending at $C_1$ and let $x$ be an element of the union of these branches. Assume $\alpha^{i}(x)=c^1_k$. Then we define $f(x)=c_{k-i-1}$. Note that this definition does not depend on the choice of $i$.

If $S_2$ and other connected components have cycles, then we map them to cycle $C^T_2$ by $f$ as in the previous case for $S_1$. Suppose that the connected component $S_2$ does not have a cycle. Then $S_2$ consists of branches. We pick up an element $b\in S_2$ and put $f(b)=c_k$. We define $f(x)=c_{k+i}$, where $x\in \alpha^{i}(b)$, $i\in \mathbb{Z}$. If there is an element $b_1\in \alpha^{-i}(\alpha^j(b))$ for which $f$ has not been defined, then we put $f(b_1)=c_{k+j-i}$.

Similarly, let $j$ be any map satisfying the conditions. Now we argue that we cannot lift it to the homomorphism $\phi$ such that $(j \circ \phi)(C_1)=f(C_1)$. We have $f(C_1)=C^T_1$, the preimage of $C^T_1$ by the map $j$ is the connected component $S_2$. As $cn(S_1)< \ cn(S_2)$, there is no map $\phi$ such that $(j \circ \phi)(C_1)=f(C_1)$.

\end{proof}

\begin{lemma}\label{mono2.2}
In a quasi-projective monounary algebra ${\cal A}=(A,\alpha)$ with a cycle does not exist an infinite branch.

The statement and the proof holds if we use the notion of factor algebra of Jakub\'ikov\'a-Studenovsk\'a's definition choosing $U=A\setminus B$ where $B$ is an infinite branch.

\end{lemma}

\begin{proof}
Suppose to the contrary that there exists an infinite branch $B$. From Lemma \ref{mono2.1} we know that it ends in a cycle $C$ and let $B=\{b_1, b_2, \ldots\}$, where $\alpha(b_1)=c_k\in C$. Let ${\cal T}=(T, \alpha_T)$ be a monounary algebra such that $T=B^T\cup C^T$ where $B^T$ is a branch in $T$ that ends at a one point cycle $C^T=\{c^T\}$ and $B^T=\{b^T_1, b^T_2, \ldots\}$.

\begin{figure}[H]
\begin{center}
\includegraphics[width= 13 cm ]{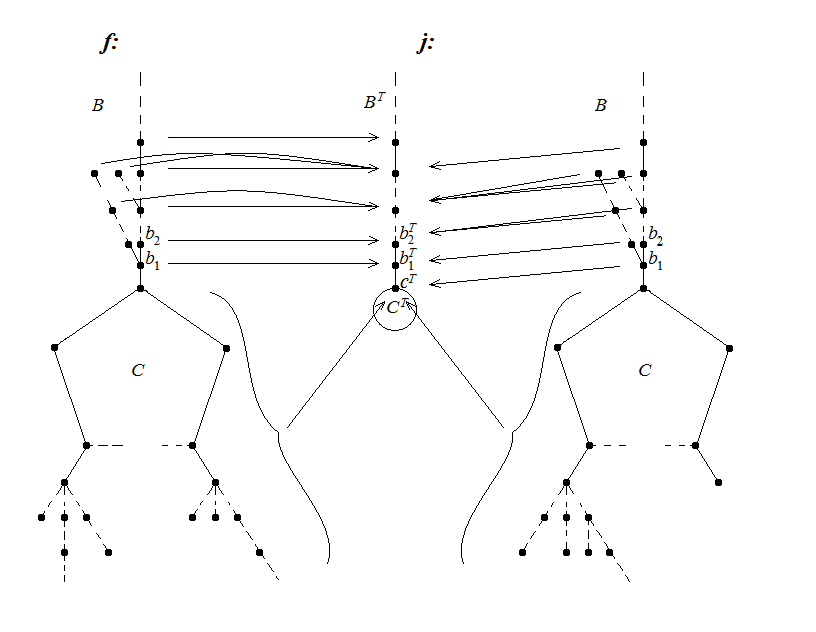}
\caption{\label{kpma1i} There is an infinite branch in a monounary algebra with a cycle }
\end{center}
\end{figure}

We define $f$ and $j$ in the following way (Fig. \ref{kpma1i}):

$$f(x)=\begin{cases} \alpha^{-i}(b^{T}_1) & if \hbox{ } x\in \alpha^{-i}(b_1), i\geq 0,
\\ c^T & otherwise \end{cases}$$

$$j(x)=\begin{cases}   \alpha^{-i}(b^{T}_{1}) & if \hbox{ } x\in \alpha^{-i}(b_{2}), i\geq 0,
\\ c^T & otherwise. \end{cases}$$

\noindent
We have $f(B)=B^T \cup C^T$ such that $f(b_1)=b^{T}_1$. Also,  $j(b_2)=b_1^{T}$. So, if there exists a map $\phi$ such that $j\circ \phi=f$, then $\phi(c_k)=b_1$  which is a contradiction.

\end{proof}

\begin{lemma}\label{mono2.3}
If a monounary algebra with a cycle is quasi-projective, then there does not exist a bunch with more than one branch.

The statement and the proof holds if we use the notion of factor algebra of Jakub\'ikov\'a-Studenovsk\'a's definition choosing $U=A\setminus (B_1\setminus B_2 \cup B_1\setminus B_2)$, where $B_1$ and $B_2$ are two different branches in a bunch.

\end{lemma}

\begin{proof} From  Lemma \ref{mono2.2} we know that there are no infinite branches. Suppose to the contrary, let $\uparrow B$ be a bunch with two different branches $B_1$ and $B_2$. From Lemma \ref{mono2.1} we know that there is a cycle $C$ that branches from $\uparrow B$ end at it. We may assume that $l_h(B_1)\geq l_h(B_2)$, $l_h(B')\leq l_h(B_1)$ and $l_h(B')\leq l_h(B_2)$ for all $B'\in \uparrow B$. Let $n+g=l_h(B_1)$, $m+g=l_h(B_2)$ and $b_g\in B_1 \cap B_2$ such that  $l_h(b_g)=l_h(B_1 \cap B_2)$, $ B_1=\{b_1, b_2\ldots, b_g, b^1_1, b^1_2,\ldots, b^1_n\}$ and $B_2\setminus B_1=\{b^2_1, b^2_2,\ldots, b^2_m\}$. Let ${\cal T}=(T, \alpha_T)$ be a monounary algebra such that $T=B^T \cup C^T$ where $B^T$ is a branch that ends at a one point cycle $C^T=\{c^{T}\}$ such that $B^T=\{b_1^{T}, b_2^{T},\ldots, b_{n}^{T}\}$. 

We define $f$ and $j$ in the following way (Fig. \ref{kpma2}):

$$f(x)=\begin{cases}   \alpha^{-i}(b^{T}_{1}) & if \hbox{ } x\in \alpha^{-i}(b^1_{1}), i\geq 0,
\\ c^T & otherwise, \end{cases}$$

$$j(x)=\begin{cases} \alpha^{-i}(b^{T}_{1}) & if \hbox{ } x\in \alpha^{-i}(b'_{1}), i\geq 0,  b_1'\in \alpha^{-1}(b_g),
 \\ c^T & otherwise. \end{cases}$$

\noindent
We have $f(B_1\cup B_2)=B^T \cup C^T$ such that $f(b^1_1)=b^{T}_1$ and $f(b^2_1)=c^T$. Also,  $j(b'_1)=b_1^{T}$ for all $b'_1\in \alpha^{-1}(b_g)$. So, if there exists a map $\phi$ such that $j\circ \phi=f$, then $\phi(b_1^1)=b^1_1$ and $\phi(b_1^2)=x$, where $x\not= b'_1$,  which is a contradiction.

\begin{figure}[H]
\begin{center}
\includegraphics[width= 13 cm ]{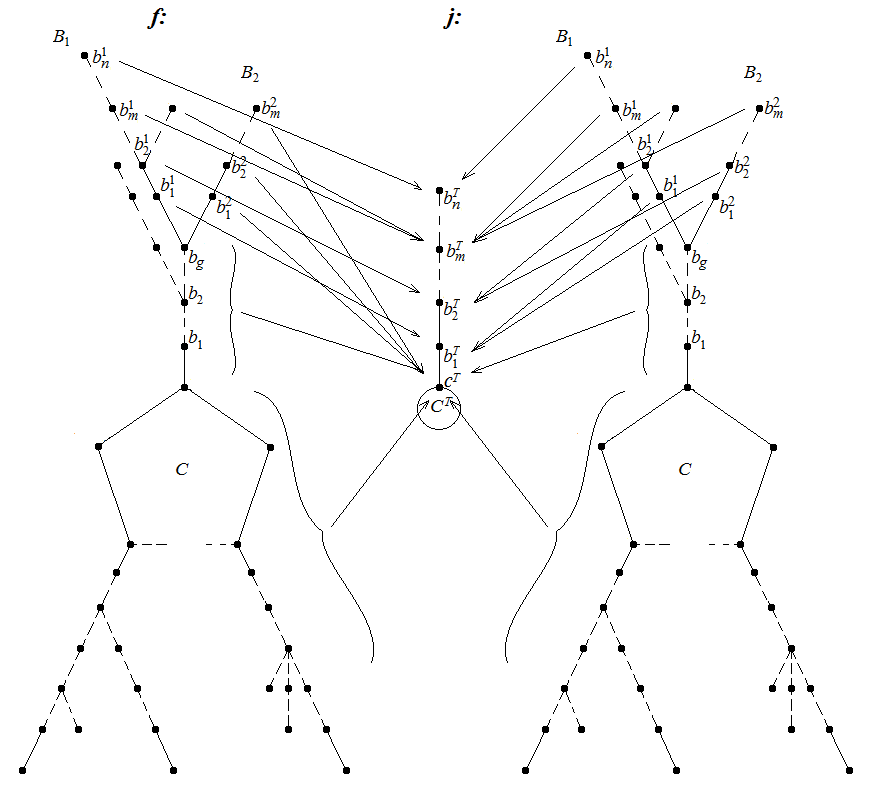}
\caption{\label{kpma2} Two branches in a bunch}
\end{center}
\end{figure}

\end{proof}

\begin{lemma}\label{mono2.4}
In a monounary algebra with a cycle all branches have the same length.

The statement and the proof holds if we use the notion of factor algebra of Jakub\'ikov\'a-Studenovsk\'a's definition choosing $U=A\setminus (B_1\cup B_2)$ where $B_1$ and $B_2$ are two branches with different lengths.

\end{lemma}

\begin{proof} Let $B_1$ and $B_2$ be branches that they end at cycles $C_1$ and $C_2$, respectively.  From Lemma \ref{mono2.1} we know that $cn(C_{1})=cn(C_{2})$. From Lemmata \ref{mono2.2} and \ref{mono2.3} we know that all branches are finite and $\uparrow B_1=B_1$ and $\uparrow B_2=B_2$. Suppose to the contrary that $l_h(B_1)\not=l_h(B_2)$. Without loss of generality we may assume that $l_h(B_1)<l_h(B_2)$.  Let $ B_1=\{ b^1_1, b^1_2,\ldots, b^1_n\}$ and $B_2=\{b^2_1, b^2_2,\ldots, b^2_m\}$. Let ${\cal T}=(T, \alpha_T)$ be a monounary algebra such that $T=B^T\cup C^T$, where $B^T$ is a branch in $T$ that ends at a one point cycle $C^T=\{c^T\}$ and $B^T=\{b_1^{T}, b_2^{T},\ldots, b_{n}^{T}\}$. We define $f$ and $j$ in the following way (Fig. \ref{kpma1}):

$$f(x)=\begin{cases} \alpha^{-i}(b^{T}_1) & if \hbox{ } x\in \alpha^{-i}(b^1_1), i\geq 0,
\\ c^T & otherwise \end{cases}$$

$$j(x)=\begin{cases}   \alpha^{-i}(b^{T}_{1}) & if \hbox{ } x\in \alpha^{-i}(b^2_{(m-n)}), i\geq 0,
\\ c^T & otherwise. \end{cases}$$

We have $f(B_1\cup C_{1})=B^T\cup C^T$ such that $f(b^1_n)=b^T_n$, the preimage of $B^T\cup C^T$ by the map $j$ is $ B_2 \cup C_{2}$ such that $j(b^2_m)=b^T_n$. Because $l_h(B_1)<l_h(B_2)$, there does not exist a map $\phi$ such that $(j \circ \phi)(B_1\cup C_{1})=f(B_1\cup C_{1})$.

\begin{figure}[H]
\begin{center}
\includegraphics[width= 12 cm ]{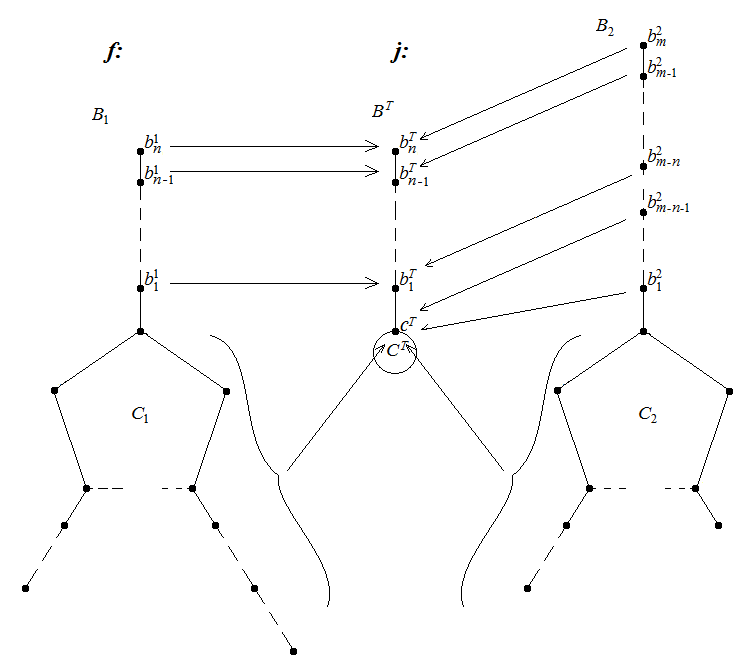}
\caption{\label{kpma1} Branches $B_1$ and $B_2$ are finite branches with different lengths }
\end{center}
\end{figure}

\end{proof}

\begin{lemma}\label{mono2.5}
Let ${\cal A}=(A,\alpha)$ be a quasi-projective monounary algebra with a connected component $S_1$ such that there are two different bunches $\uparrow B' $ and $\uparrow B''$ in $S_1$ which end at a cycle $C$ with $cn(C)=m$ in elements $c_k$ and $c_l$, $k<l$, respectively, where $\alpha^{l-k}(c_k)=c_l$, $\alpha^{m-l+k}(c_l)=c_k$. Then $m$ is even and $l-k= m/2$.

The statement and the proof holds if we use the notion of factor algebra of Jakub\'ikov\'a-Studenovsk\'a's definition choosing $U=A\setminus (B_1\cup B_2)$, where $B_1$ and $B_2$ are the maximal branches in $\uparrow B' $ and $\uparrow B''$, respectively.
\end{lemma}

\begin{proof}
Suppose to the contrary that $m$ is odd or $m \not=2l-2k$. Let $B_1 \in \uparrow B' $ and $B_2\in \uparrow B''$. From Lemmata \ref{mono2.2}, \ref{mono2.3} and \ref{mono2.4} we know that branches are finite, $\uparrow B'=B_1$, $\uparrow B''=B_2$, $l_h(B_1)=l_h(B_2)$ and all branches have the same length. Let $ B_1=\{ b^1_1, b^1_2,\ldots, b^1_n\}$ and $B_2=\{b^2_1, b^2_2,\ldots, b^2_n\}$. 

\begin{figure}[H]
\begin{center}
\includegraphics[width= 13 cm ]{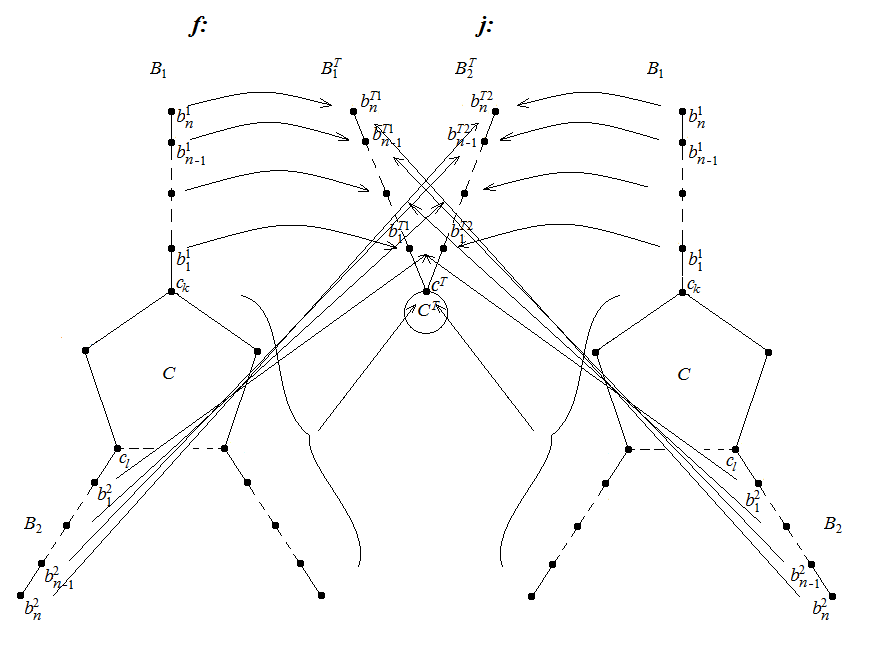}
\caption{\label{kpma4} $B_1$ and $B_2$ end at elements $c_k$ and $c_l$, respectively  }
\end{center}
\end{figure}

Let ${\cal T}=(T, \alpha_T)$ be a monounary algebra such that $T=B^T_1 \cup B^T_2 \cup C^T$ where $B^T_1$ and $B^T_2$ are branches that ends at a one point cycle $C^T=\{c^{T}\}$, $B_1^T=\{b_1^{T1}, b_2^{T1},\ldots, b_{n}^{T1}\}$ and $B_2^T=\{b_1^{T2}, b_2^{T2},\ldots, b_{n}^{T2}\}$. We define $f$ and $j$ in the following way  (Fig. \ref{kpma4}):

$$f(x)=\begin{cases}  \alpha^{-i}(b^{T1}_{1}) & if \hbox{ } x\in \alpha^{-i}(b^1_1), i\geq 0,
\\ \alpha^{-i}(b^{T2}_{1}) & if \hbox{ } x\in \alpha^{-i}(b^2_1), i\geq 0,
 \\ c^T & otherwise, \end{cases}$$

$$j(x)=\begin{cases}  \alpha^{-i}(b^{T2}_{1}) & if \hbox{ } x\in \alpha^{-i}(b^1_1), i\geq 0,
\\ \alpha^{-i}(b^{T1}_{1}) & if \hbox{ } x\in \alpha^{-i}(b^2_1), i\geq 0,
 \\ c^T & otherwise. \end{cases}$$

\noindent
We have $f(B_ 1\cup B_2 \cup C)=B^T_1 \cup B^T_2 \cup C^T_{B}$ such that $f(c_{k})=c^T$, the preimage of $B^T_1 \cup B^T_2 \cup  C^T_{B}$ by the map $j$ is $B_1 \cup B_2 \cup C$ such that $f(c_{l})=c^T$. We have to map $B_1$ to $B_2$ and $B_2$ to $B_1$ with the cycle $C$ by $\phi$, but because $2l-2k\not= m$, there does not exist a map $\phi$ such that $(j\circ \phi)(C)=f(C)$.

\end{proof}

\begin{lemma}\label{mono2.6}
If a monounary algebra with a cycle is  quasi-projective, then there are at most two disjoint branches in a connected component or in every connected component there is at most one branch.

The statement and the proof holds if we use the notion of factor algebra of Jakub\'ikov\'a-Studenovsk\'a's definition choosing $U=A\setminus (B_1\cup B_2\cup B_3)$ where $B_1$, $B_2$ and $B_3$ are branches such that at least two are in the same connected component.
\end{lemma}

\begin{proof}

Let  $B_1$, $B_2$ and $B_3$ be three branches.   From Lemmata \ref{mono2.2}, \ref{mono2.3} and \ref{mono2.4} we know that branches are finite, every bunch has just one branch, $l_h(B_1)=l_h(B_2)=l_h(B_3)$ and all branches have the same length. From Lemma \ref{mono2.5} we know that there cannot be three branches in a connected component. Suppose to the contrary and let  branches $B_1$ and $B_2$ be in a same connected component. Let $C_1$ and $C_2$ be two cycles such that $B_1$ and $B_2$ end at the cycle $C_1$ and $B_3$ ends at the cycle $C_2$. From Lemma \ref{mono2.1} we  know that $cn(C_{1})=cn(C_{2})$. Also, we know that the ending points of branches $B_1$ and $B_2$ are on the same distance from each other.

Let $B_1=\{b^1_1,\ldots,b^1_n\}$, $B_2=\{b^2_1,\ldots,b^2_n\}$ and $B_3=\{b^3_1,\ldots,b^3_n\}$. Let ${\cal T}=(T, \alpha_T)$ be a monounary algebra such that $T=B^T_1 \cup B^T_2 \cup C^T$ where $B^T_1$ and $B^T_2$ are branches that ends at a one point cycle $C^T=\{c^{T}\}$, $B^T_1=\{b^{T1}_1,b^{T1}_2,\ldots, b^{T1}_n\}$, $B^T_2=\{b^{T2}_1,b^{T2}_2,\ldots,b^{T2}_n\}$. We define $f$ and $j$ in the following way  (Fig. \ref{kpma3eV}):

$$f(x)=\begin{cases}  \alpha^{-i}(b^{T1}_{1}) & if \hbox{ } x\in \alpha^{-i}(b^1_1) , i\geq 0,
\\ \alpha^{-i}(b^{T2}_{1}) & if \hbox{ } x\in \alpha^{-i}(b^2_1), i\geq 0,
 \\ c^T & otherwise, \end{cases}$$

$$j(x)=\begin{cases}  \alpha^{-i}(b^{T2}_{1}) & if \hbox{ } x\in \alpha^{-i}(b^1_1) , i\geq 0,
\\ \alpha^{-i}(b^{T1}_{1}) & if \hbox{ } x\in \alpha^{-i}(b^3_{1}), i\geq 0,
 \\ c^T & otherwise. \end{cases}$$

\noindent
We have $f(B_1 \cup B_2 \cup C_{1})=B^T_1 \cup B^T_2 \cup C^T_{B}$, the preimage of $B^T_1 \cup B^T_2 \cup  C^T$ by the map $j$ is in two different connected components, so there does not exist a map $\phi$ such that $(j\circ \phi)(B_1 \cup B_2 \cup C_{1})=f(B_1 \cup B_2 \cup C_{1})$.

\begin{figure}[H]
\begin{center}
\includegraphics[width= 13 cm ]{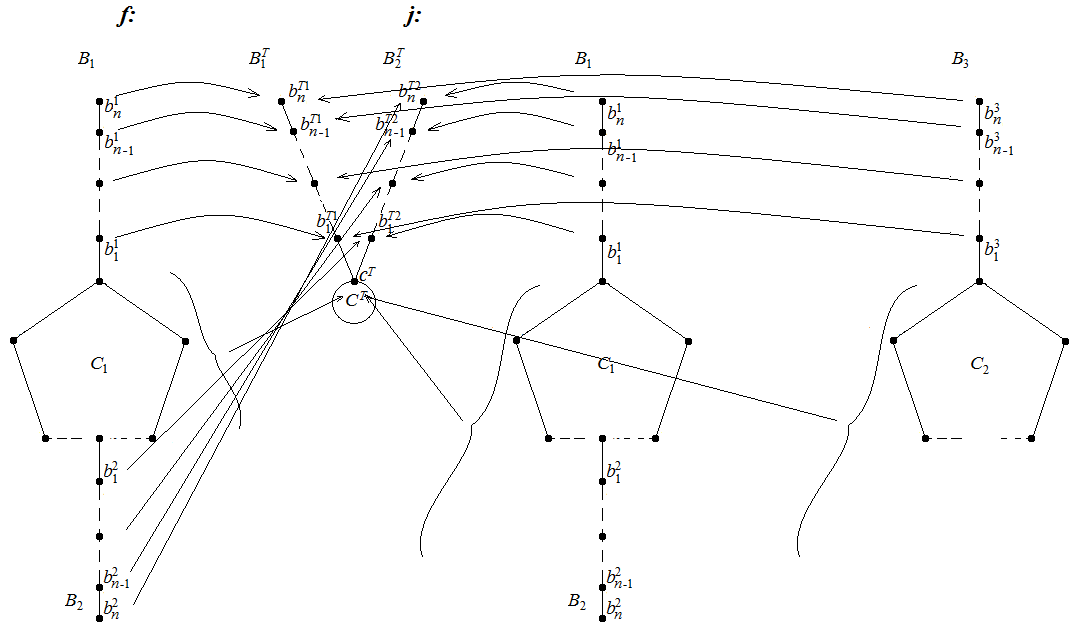}
\caption{\label{kpma3eV} Proof of Lemma \ref{mono2.6}}
\end{center}
\end{figure}

\end{proof}

\begin{lemma}\label{mono2.21} Let ${\cal A}=(A,\alpha)$ be a monounary algebra with a following:

\begin{enumerate}

\item there are cycles in all connected components,

\item cycles have the same length  in all connected components,

\item in all bunches there is one branch,

\item all branches have the same length,

\item 
\begin{enumerate}
\item[(i)] in every connected component there is at most one branch or 
\item[(ii)] there are two disjoint branches in a connected component and their cyclic element at which they end is on the same distance from each other.
\end{enumerate}

\end{enumerate}

Then ${\cal A}=(A,\alpha)$ is quasi-projective.

The statement holds if we use the notion of factor algebra of Jakub\'ikov\'a-Studenovsk\'a's definition.

\begin{figure}[H]
\begin{center}
\includegraphics[width= 14cm ]{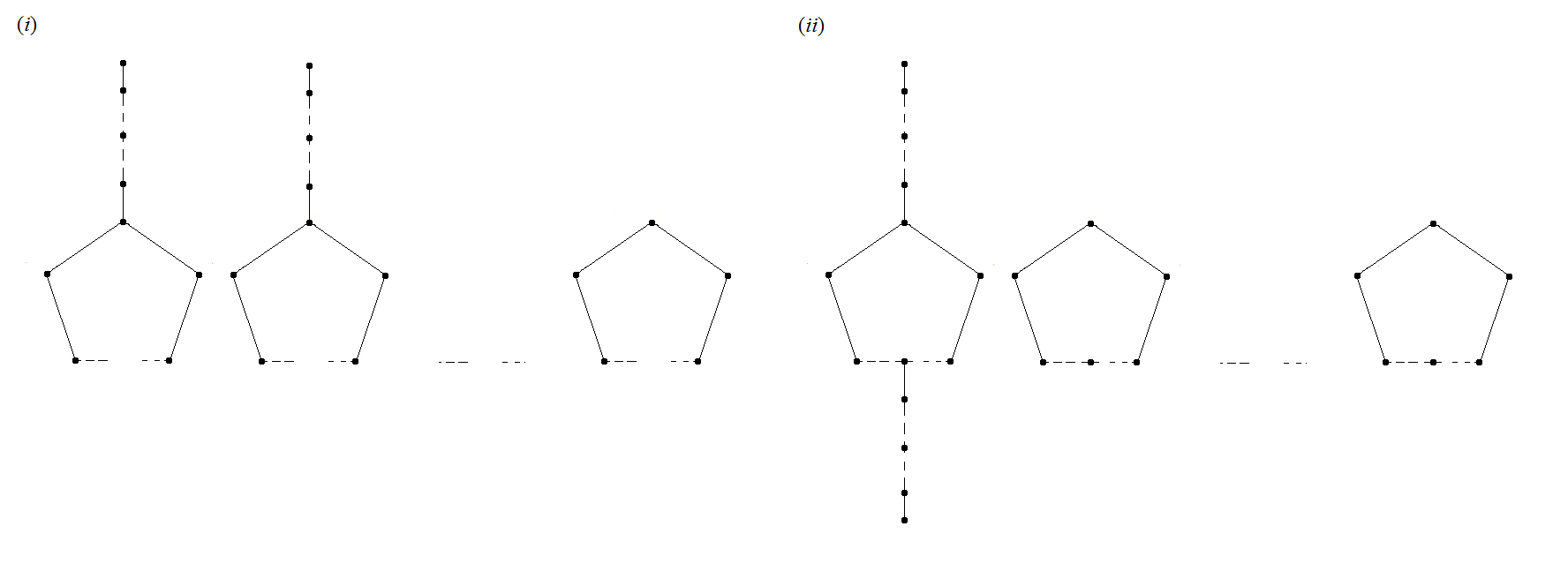}
\caption{\label{Tetel41} Lemma \ref{mono2.21}}
\end{center}
\end{figure}
\end{lemma}

\begin{proof}

Let ${\cal A}=(A,\alpha)$ be a monounary algebra where the previous conditions hold. Let ${\cal T}=(T,\alpha_T)$ be an arbitrary monounary algebra, $f:{\cal A}\rightarrow {\cal T}$ be a homomorphism and $j:{\cal A}\rightarrow {\cal T}$ be an epimorphisms.

Suppose ${\cal A}=(A,\alpha)$ has a cycle $C_1$. Then $f(C_1)\subseteq j(A)$, because $j$ is an epimorphism. $f(C_1)$ is a cycle, where $cn(f(C_1))|\ cn(C_1)$, with or without branches that end at it. The preimage of $f(C_1)$ by the map $j$ can be branches which end at cycles or just simply cycles. Infinite branches cannot occur, because then it contradicts the assumption. It follows from assumptions \ref{item11} and \ref{item12} that there is at least one cycle $C_2$ such that $j(C_2)=f(C_1)$ and  $cn(C_1)=\ cn(C_2)$. If there does not exist a branch $B_1$ which ends at the cycle $C_1$, then we define a map for the elements from $ C_1 $ to $ C_2 $ in the following way: $\phi(c^1_{i})=j^{-1}(f(c^1_{i}))$, where $j^{-1}(f(c^1_{i})) \in  C_2$.

 Now, suppose there is a branch $B_1$ which ends at the cycle $C_1$. We know that $\uparrow B_1 =B_1$ from assumption \ref{item13}. Suppose that there is no other bunch in this connected component. $f({  B_1 })\subseteq j(A)$, because $j$ is an epimorphism. Without loss of generality, let $c_1\in B_1\cap C_1$.  If $f(  B_1 )=f(C_1)$, then we find $f(c_1)$ in the $Im(j)$ and define the following mapping for the elements from $  B_1 $ to $C_2$ and from $C_1$ to $C_2$:  $\phi(b_i)=j^{-1}(f(b_i))$, where $j^{-1}(f(b_i))\in C_2$ and $\phi(c^1_{i})=j^{-1}(f(c^1_{i}))$, where $j^{-1}(f(c^1_{i})) \in  C_2$. If $f( B_1 ) \cap f(C_1) \not=\emptyset$, then for the part of $ B_1 $ which is mapped to $f(C_1)$ by the map $f$ we make the same mapping for $\phi$ as previously. So, without loss of generality, we may assume that $f( B_1 ) \cap f(C_1) =\emptyset$. There is a bunch  $ \uparrow  B_2 $ which ends at the cycle $C_2$ or there is another cycle $C_3$ with the same length as $C_1$ and a bunch $ \uparrow  B_3 $ which ends at the cycle $C_3$ unless the map $j$ is not epimorphism.

Without loss of generality, assume that that there is a bunch $\uparrow  B_2 $ which ends in the cycle $C_2$ and which is mapped to the bunch $f( B_1 )$ by the map $j$. From assumption \ref{item13} we know that $\uparrow B_2=B_2$. We define a map for the elements from $ B_1 $ to $ B_2 $ and from $C_1$ to $C_2$ in the following way: $\phi(b^1_{i})=j^{-1}(f(b^1_{i}))$, where $j^{-1}(f(b^1_{i})) \in  B_2$ $\phi(c^1_{i})=j^{-1}(f(c^1_{i}))$, where $j^{-1}(f(c^1_{i})) \in  C_2$. From assumption \ref{item14} we know that all branches have the same length. So, it cannot happen that we have to map leaf from $B_1$ to leaf from $B_2$ such that the branch $B_2$ is longer than the branch $B_1$.

Now, suppose that there is an another bunch $\uparrow  B'_1 $ which ends at the cycle $C_1$. From assumption \ref{item13} and \ref{item15} we know that then $\uparrow  B_1 $ and $\uparrow  B'_1 $ are just branches $B_1$ and $B'_1$, they are on the same distance and there are not other branches in $\cal A$. If $f(B_1 )=f(C_1)$ or $f(B'_1)=f(C_1)$, we handle it as if there were just one branch, so suppose $f(B_1 )\not=f(C_1)$ and $f(B'_1)\not=f(C_1)$. Without loss of generality, we suppose $f(B_1) \cap f(C_1) =\emptyset$ and $f( B'_1) \cap f(C_1) =\emptyset$. From $f(B_1) \cap f(C_1) =\emptyset$ and $f(B'_1) \cap f(C_1) =\emptyset$ and because $j$ is an epimorphism, we have that there are bunches $\uparrow  B_2 $ and $\uparrow  B'_2 $ which end at a cycle $C_2$ such that $j(\uparrow  B_2 )=f(B_1)$, $j(\uparrow  B'_2 )=f(B'_1)$, $j(C_2)=f(C_1)$, where cn~$(C_1)=$\ cn~$(C_2)$. Because of the assumption that there are not any other branches, we have that $\uparrow B_2= B_2 $ and $\uparrow B'_2= B'_2 $, so $\{ B_1 , B'_1\}=\{B_2 , B'_2\}$.  Let $ B_1= B_2$ and $ B'_1= B'_2 $. We define a map for the elements from $ B_1 $ and $ B'_1 $ to $ B_2$ and $B'_2 $, respectively,  in the following way: $\phi(b^1_{i})=j^{-1}(f(b^1_{i}))$ and $\phi(b^{'1}_{i})=j^{-1}(f(b^{'1}_{i}))$, where $j^{-1}(f(b^1_{i})) \in B_2$ and $j^{-1}(f(b'^1_{i}))\in B'_2$, respectively. Also, the mapping from $C_1$ to $C_2$ is the following: $\phi(c^1_{i})=j^{-1}(f(c^1_{i}))$, where $j^{-1}(f(c^1_{i})) \in  C_2$. The ending points of the branches $B_1$ and $ B'_1$ are on the same distance from each other, so it cannot happen we cannot map the cycle $C_1$ to itself. 

If we use the notion of factor algebra of Jakub\'ikov\'a-Studenovsk\'a's definition, then the proof is analogous.

\end{proof}

\section{Monounary algebras without cycles}

In this subsection we characterise quasi-projective monounary alegbras without cycles.

\begin{lemma}\label{mono2.91}

If a monounary algebra without a cycle is quasi-projective, then there does not exist two branches in a connected component.

The statement and the proof holds if we use the notion of factor algebra of Jakub\'ikov\'a-Studenovsk\'a's definition choosing $U=A\setminus (B_1\setminus B_2 \cup B_1\setminus B_2)$, where $B_1$ and $B_2$ are two different branches in a connected component.

\end{lemma}

\begin{proof} Suppose to the contrary, let $S_1$ be a connected component with two different branches $B'$ and $B''$. Let $b_g\in B' \cap B''$ such that $l_h(b_g)=l_h(B' \cap B'')$. Let $B_1$ and $B_2$ be branches in $\uparrow b_g$ such that $l^g_d(B)\leq l^g_d(B_1)$ and $l^g_d(B)\leq l^g_d(B_2)$ for all $B\in \uparrow b_g$. We may assume that $l^g_d(B_1)\geq l^g_d(B_2)$. Let $b^1_1 \in B_1$ and $b^2_1\in B_2$ such that $\alpha(b^1_1)=\alpha(b^2_1)=b_g$. Let ${\cal T}=(T, \alpha_T)$ be a monounary algebra such that $T=B^T \cup C^T$ where $B^T$ is a branch that ends at a one point cycle $C^T=\{c^{T}\}$ such that $l_d(B^T)=l_d(B_1)$, $b^T_1\in B^T$ and $\alpha(b^T_1)=c^T$. We define $f$ and $j$ in the following way (Fig. \ref{kpma2}):

$$f(x)=\begin{cases}   \alpha^{-i}(b^{T}_{1}) & if \hbox{ } x\in \alpha^{-i}(b^1_{1}), i\geq 0,
\\ c^T & otherwise, \end{cases}$$

$$j(x)=\begin{cases} \alpha^{-i}(b^{T}_{1}) & if \hbox{ } x\in \alpha^{-i}(b'_{1}), i\geq 0, b_1'\in \alpha^{-1}(b_g),
 \\ c^T & otherwise. \end{cases}$$

\noindent
We have $f(B_1\cup B_2)=B^T \cup C^T$ such that $f(b^1_1)=b^{T}_1$ and $f(b^2_1)=c^T$. Also,  $j(b'_1)=b_1^{T}$ for all $b'_1\in \alpha^{-1}(b_g)$. So, if there exists a map $\phi$ such that $j\circ \phi=f$, then $\phi(b_1^1)=b^1_1$ and $\phi(b_1^2)=x$, where $x\not= b'_1$,  which is a contradiction.

\end{proof}

\begin{figure}[H]
\begin{center}
\includegraphics[width= 15 cm ]{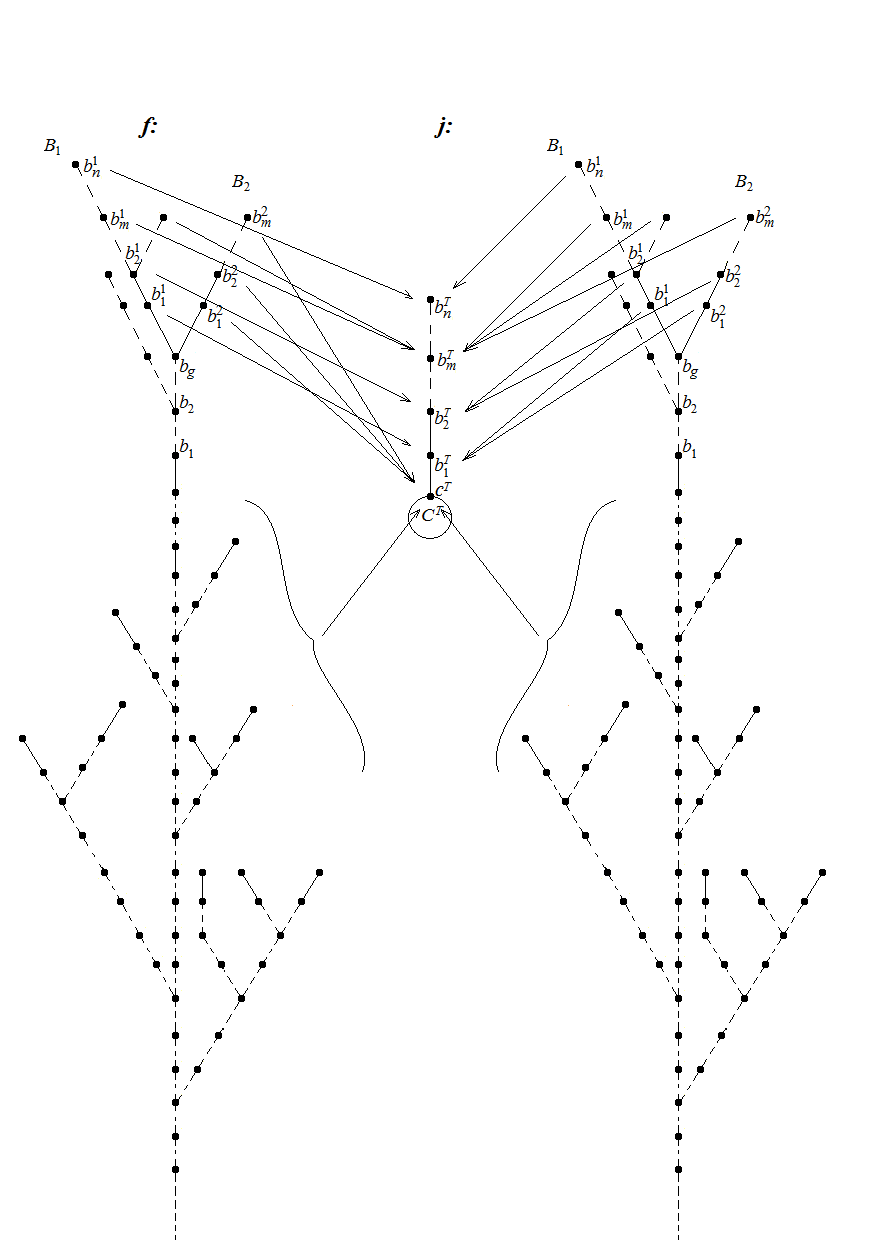}
\caption{\label{kpma6} Proof of Lemma \ref{mono2.91}}
\end{center}
\end{figure}

\begin{theorem}\label{mono2.31}
Let ${\cal A}=(A,\alpha)$ be a monounary algebra such that in every connected component there is just one branch. Then it is quasi-projective.

The statement holds if we use the notion of factor algebra of Jakub\'ikov\'a-Studenovsk\'a's definition.

\end{theorem}

\begin{proof} Let ${\cal A}=(A,\alpha)$ be a monounary algebra such that in every connected component there is just one branch. Let ${\cal T}=(T,\alpha_T)$ be an arbitrary monounary algebra, $f:{\cal A}\rightarrow {\cal T}$ be a homomorphism and $j:{\cal A}\rightarrow {\cal T}$ be an epimorphisms. Let $B$ be a branch. Then $f(B)\subseteq j(A)$, because $j$ is an epimorphism. $f(B)$ can be a cycle or a branch or both. Let $B'$ be a connected component such that $j(B')=f(B)$. We define the map for the elements from $ B $ to $ B' $ in the following way: $\phi(b_{i})=j^{-1}(f(b_{i}))$, where $j^{-1}(f(b^1_{i})) \in  B'$. It can be seen trivially that $\phi$ is a homomorphism. 

If we use the notion of factor algebra of Jakub\'ikov\'a-Studenovsk\'a's definition, then the proof is analogous.

\end{proof}

\noindent
{\bf Acknowledgement.} The author would like to express gratitude to  Csaba Szab\'o  and G\'abor Somlai  for their careful reading, comments and remarks. They highly improved the presentation of this article.

\bibliographystyle{plain}

\vskip5mm
\noindent
{\sc \' Eva Jung\'abel, E\"otv\"os Lor\'and University, P\'azm\'any P\'eter s\'et\'any 1/C, 1117 Budapest, Hungary}

\noindent
{\it E-mail address:} {\tt evajungabel@student.elte.hu}

 \end{document}